\documentclass{elsart}

\usepackage[numbers]{natbib}         % Bibliography
\usepackage[english]{babel}
\usepackage{amsfonts}

\newcommand{\p}{\mathbb{P}}          
\newcommand{\E}{\mathbb{E}}          
 
\newcommand{\A}{\mathcal{A}} 
\newcommand{\n}{\mathcal{N}} 

\newcommand{\one}{\bf 1}

\begin{document}

\begin{frontmatter}

\title{Poisson approximation for search of rare words in DNA sequences}

\author[Evry]{Nicolas Vergne\corauthref{cor}}
\corauth[cor]{Corresponding author.}
\ead{nicolas.vergne@genopole.cnrs.fr}
\author[Sao]{Miguel Abadi}
\ead{abadi@ime.unicamp.br}

\address[Evry]{Universit\'{e} d'Evry Val d'Essonne, D\'{e}partement Math\'{e}matiques,
 Laboratoire Statistique et G\'{e}nome,
91 000 Evry, France}
\address[Sao]{Universidade de Campinas, Brasil}

\begin{abstract}
Using recent results on the occurrence times of a string of 
symbols in a stochastic process with mixing properties, 
we present a new method for the search of rare words in
biological sequences generally modelled by a Markov chain. We obtain a bound on the error between 
the distribution of the number of occurrences of a word in a sequence (under a Markov model) 
and its Poisson approximation. 
A global bound is already given by a Chen-Stein method. Our approach, 
the $\psi$-mixing method, gives local bounds. 
Since we only need  the error in the tails of 
distribution, the global uniform bound of Chen-Stein is too large and it is a better
way to consider local bounds.
We search for two thresholds on the number of occurrences from which we can regard the studied word as an
over-represented or an under-represented one. 
A biological role is suggested for these over-
or under-represented words. 
Our method gives such thresholds for a panel of words much broader than the Chen-Stein method.
Comparing the methods, we observe a better accuracy for the $\psi$-mixing method for the bound
of the tails of distribution. We also present the software 
\texttt{PANOW}\footnote{available at \texttt{http://stat.genopole.cnrs.fr/software/panowdir/}}
dedicated to the computation of the error term and the thresholds for a studied word.
\end{abstract}

\begin{keyword}
  Poisson approximation\sep
  Chen-Stein method\sep 
  mixing processes\sep 
  Markov chains\sep
  rarewords\sep
  DNA sequences
\end{keyword}

\end{frontmatter}

%INTRODUCTION

\section{Introduction}

Modelling DNA sequences with stochastic models and developing statistical methods to analyse the
enormous set of data that results from the multiple projects of DNA sequencing are challenging
questions for statisticians and biologists.
Many DNA sequence analysis are based on the distribution of the occurrences of patterns having some
special biological function. The most popular model in this domain is the Markov chain model that
gives a description of the local behaviour of the sequence (see
\citet{MC, MCA, hexa, Extend}).
An important problem is to determine the statistical significance of a word frequency in a DNA
sequence. 
\citet{Nicod} discuss about this relevance of finding over- or under-represented words.
The naive idea is the following: a word may have a significant low frequency in a DNA sequence
because it disrupts replication or gene expression, whereas a significantly frequent word may have a
fundamental activity with regard to genome stability. Well-known examples of words with exceptional
frequencies in DNA sequences are biological palindromes corresponding to restriction sites
avoided for instance in \textit{E. coli} (\citet{Karlin2}), the Cross-over Hotspot Instigator
sites in several bacteria, again in \textit{E. coli} for example (\citet{ecoli,ElK}), and uptake
sequences (\citet{uptake}) or polyadenylation signals (\citet{Van_H}).

The exact distribution of the number of a word occurrences under the
Markovian model is known and some softwares are available (\citet{jjd, MR}) but, because of numerical
complexity, they are often used 
to compute  
expectation and variance of a given count (and thus use, in fact, Gaussian approximations for the
distribution). 
In fact these methods are not efficient for long sequences or if the Markov model order is larger
than $2$ or $3$.
For such cases, several approximations are possible:
Gaussian approximations (\citet{Bernard2}), Binomial (or Poisson) approximations
(\citet{Van_H2,Poisson}), 
compound Poisson approximations (\citet{ss2}), or large deviations approach (\citet{GregSpatt}).
In this paper we only focus on the Poisson approximation. 
We approximate $\p(N(A)=k)$ by $\exp(-t\p(A))[t\p(A)]^k(k!)^{-1}$
where $\p(N(A)=k)$ is the stationary probability under the Markov model that the number of
occurrences $N(A)$ of  word $A$ is equal to $k$, $\p(A)$ is the probability that word $A$
occurs at a given position, and $t$ is the length of the sequence. 
Intuitively, a binomial distribution could be used to approximate the distribution of
occurrences of a particular word. Length $t$ of the sequence is large, $\p(A)$ is small and $t\p(A)$
is almost constant. Thus, we
use the more numerically convenient Poisson approximation.
Our aim is to bound the error between the distribution of the number of occurrences of word $A$
and its Poisson approximation. 
In \citet{ss2}, the authors prove an upper bound for 
a compound Poisson approximation. They use a Chen-Stein method, which is the
usual method in this purpose. This method has been developed by Chen on Poisson
approximations (\citet{Chen}) after a work of Stein on normal approximations (\citet{Stein}).
Its principle is to bound the difference between the two distributions in
total variation distance for all subsets of the definition domain. 
Since we are interested in under- or over-represented words, we are only
interested in this difference for the tails of the
distributions. Then, the uniform bound given by the Chen-Stein method is too large for our purpose.
We present here a new method, based on the property 
of mixing processes. Our method has the useful particularity to give a bound on the error at each
point of the distribution. More precisely, it offers an error term $\epsilon$, for the number of
occurrences $k$, of word $A$:  
\[
\left|\p(N(A)=k)-\frac{{e}^{-t\p(A)}{(t\p(A))}^k}{k!}\right|
 \leq \epsilon(A,k).
\]
Moreover, $\epsilon(A,k)$ decays factorially fast with respect to $k$.

\citet{abadi2, abadi3} presents lower and upper bounds for the exponential approximation of the
first occurrence time of a rare event, also called \textit{hitting time}, in a stationary stochastic
process on a finite alphabet with $\alpha$- or $\phi$-mixing property. \citet{abadi1} describe the
statistics of \textit{return times} of a string of symbols in such a process. In \citet{abadi4}, the
authors prove a Poisson approximation for the distribution of occurrence times of a string of
symbols in a $\phi$-mixing process. The first
part of our work is to determine some constants not explicitly computed in the results of the
above mentioned articles but necessary for the proof of our theorem. 
Our work is complementary to all these articles, in the sense
that it relies on them for preliminary results and it adapts them to
$\psi$-mixing processes. Since Markov chains 
are mixing processes, all these results established for mixing processes also apply to Markov
chains which model biological sequences.

This paper is organised in the following way. In section \ref{CS}, we introduce the Chen-Stein
method. In section \ref{MA}, we define a $\psi$-mixing process and state some preliminary notations,
mostly on the properties of a word. We also present in this section the principal result of our
work: the Poisson approximation (Theorem \ref{bigt2}). 
In section \ref{Prel}, we state preliminary results. Mainly, we recall results of \citet{abadi3},
computing all the necessary constants and we present lemmas and propositions necessary for the
proof of Theorem \ref{bigt2}.
In section \ref{Fin}, we establish the proof of our main result: Theorem \ref{bigt2} on Poisson
approximation. Using $\psi$-mixing properties and preliminary results, we prove an upper bound for
the difference between the exact distribution of the number of occurrence of word $A$ and the
Poisson distribution of parameter $t\p(A)$.
Section \ref{BIO} is dedicated to numerical results. For the search of
over-represented words, we compare our method to Chen-Stein method
on both synthetic and biological data. In this section, we also present results obtained by a
similar method, the $\phi$-mixing method.
We end the paper presenting some examples of biological applications, and some conclusions and
perspectives of future works.

%CHEN-STEIN

\section{The Chen-Stein method}\label{CS}

\subsection{Total variation distance}

\begin{defn} For any two random variables $X$ and $Y$ with values in the same discrete space $E$,
the total variation distance between their probability distributions is defined by
\[d_{\mathrm{TV}}(\mathcal{L}(X), \mathcal{L}(Y))=\frac{1}{2}\sum_{i \in E}
\left|\p(X=i)-\p(Y=i)\right|.\]
\end{defn}
We remark that for any subset $S$ of $E$
\[\left|\p(X \in S)-\p(Y \in S)\right| \leq d_{\mathrm{TV}}(\mathcal{L}(X), 
\mathcal{L}(Y)).\]

\subsection{The Chen-Stein method}

The Chen-Stein method is used to bound the error between the distribution of the number
of occurrences of a word $A$ in a sequence $X$ and the Poisson distribution with parameter $t\p(A)$ 
where $t$
is the length of the sequence and $\p(A)$ the stationary measure of $A$.
The Chen-Stein method for Poisson approximation has been developed by \citet{Chen}; a friendly
exposition is in \citet{Ar1} and a description with many examples can be found in \citet{Ar2} and
\citet{Bar}. We will use Theorem $1$ in \citet{Ar2} with an improved bound by \citet{Bar} (Theorem
$1$.A and Theorem $10$.A).

First, we will fix a few notations.
Let $\mathcal{A}$ be a finite set (for example, in the DNA case $\mathcal{A}=\{a,c,g,t\}$).
Put $\Omega=\mathcal{A}^{\Zset}$.
For each
 $x={\left(x_m\right)}_{m\in{\Zset}}\in{\Omega}$, we denote by $X_m$ the $m$-th coordinate
of the sequence $x$: $X_m(x)=x_m$. We denote by 
$T: \Omega\rightarrow \Omega$ the one-step-left shift operator: so we will have
${\left(T(x)\right)}_m=x_{m+1}$. 
We denote by $\mathcal{F}$ the $\sigma$-algebra over $\Omega$ generated by strings and by
 $\mathcal{F}_I$ the  $\sigma$-algebra generated by strings with coordinates in $I$
with $I\subseteq{\Zset}$. We consider an invariant probability measure $\p$ 
over $\mathcal{F}$. 
Consider a stationary Markov chain $X=\left(X_i\right)_{i \in \Zset}$ on the finite alphabet $\mathcal{A}$.
Let us fix a word $A= (a_1, \dots, a_n)$. For $i\in\{1,2,\cdots,t-n+1\}$,
let $Y_i$ be the following random variable 
\begin{eqnarray}
Y_i=Y_i(A) & = & \one \{\mbox{word $A$ appears at position $i$ in the sequence}\} \nonumber \\ 
           & = & \one \{ (X_i, \dots, X_{i+n-1})=(a_1, \dots, a_n)\}, \nonumber
\end{eqnarray}
where $\one\{F\}$ denotes the indicator function of set $F$.
We put $Y=\sum_{i=1}^{t-n+1}Y_i$, the random variable corresponding to the number of
occurrences of a word, $\E(Y_i)=m_i$ and 
$\sum_{i=1}^{t-n+1}m_i=m$. 
Then, $\E(Y)=m$.
Let $Z$ be a Poisson random variable with parameter $m$: $Z\sim\mathcal{P}(m)$. 
For each $i$, we arbitrarily define a set $V(i)\subset \{1,2,\cdots,t-n+1\}$ containing the point $i$.
The set  $V(i)$ will play the role of a neighbourhood of $i$.

\begin{thm}[\citet{Ar2, Bar}]\label{ChenStein}
Let I be an index set. For each $i \in I$, let $Y_i$ be a Bernoulli random variable with
$p_i=\p(Y_i=1)>0$. Suppose that, for each $i \in I$, we have chosen $V(i) \subset I$ with $i \in V(i)$.
Let $Z_i, i \in I$, be independent Poisson variables with mean $p_i$. The total variation distance
between the dependent Bernoulli process $\underline{Y}=\{Y_i, i \in I\}$ and the Poisson process
$\underline{Z} =\{Z_i, i\in I\}$ satisfies
\[d_{\mathrm{TV}}(\mathcal{L}(\underline{Y}),\mathcal{L}(\underline{Z}))\leq b_1+b_2+b_3\]
where
\[b_1=\sum_i\sum_{j\in{V(i)}}\E(Y_i)\E(Y_j),\]
\[b_2=\sum_i\sum_{j\in{V(i)},  j\neq{i}}\E(Y_iY_j),\]
\[b_3=\sum_i\E\left|\E(Y_i-p_i  |  Y_j,j\notin{V(i)})\right|.\]
Moreover, if $W=\sum_{i \in I} Y_i$ and $\lambda=\sum_{i \in I} p_i< \infty$, then 
\[d_{TV}(\mathcal{L}(W), \mathcal{P}(\lambda))\leq\frac{1-e^{-\lambda}}{\lambda}(b_1+b_2)+ \min{\left(1,
  \sqrt{\frac{2}{\lambda e}}\right)}b_3.\]
\end{thm} 

We think of $V(i)$ as a neighbourhood of strong dependence of $Y_i$.
Intuitively, $b_1$ describes the contribution related to the size of the neighbourhood and the
weights of the random variables in that neighbourhood; if all $Y_i$ had the same probability of
success, then $b_1$ would be directly proportional to the neighbourhood size. The term $b_2$ accounts
for the strength of the dependence inside the neighbourhood; as it depends on the second moments, it
can be viewed as a ``second order interaction'' term. Finally, $b_3$ is related to the strength of
dependence of $Y_i$ with random variables outside its neighbourhood. In particular, 
note that $b_3=0$ if $Y_i$ is independent of $\left\{Y_j | j \notin V(i)\right\}$. 

One consequence of this theorem is that for any indicator function of an event, i.e. for any measurable
functional $h$ from $\Omega$ to $[0,1]$, there is an error bound of the form $|\E h(\underline{Y}) -
\E h(\underline{Z})| \leq d_{TV}(\mathcal{L}(\underline{Y}),\mathcal{L}(\underline{Z}))$. Thus, if
$S(\underline{Y})$ is a test statistic then, for all $t \in \Rset$,
\[\p(S(\underline{Y}) \geq t) -\p(S(\underline{Z}) \geq t) \leq b_1+b_2+b_3,\]
which can be used to construct confidence intervals and to find p-values for tests based on this
statistic.

%LOI DE POISSON
\section{Preliminary notations and Poisson Approximation}\label{MA}

\subsection{Preliminary notations}\label{prelim}

We focus on Markov processes in our biological applications (see \ref{BIO}) but the theorem given in the 
following subsection is established for more general mixing processes: the so called 
$\psi$-mixing processes. 
\begin{defn}
Let $\psi={(\psi(\ell))}_{\ell\geq{0}}$ be a sequence of real numbers decreasing
to zero. We say that 
 ${(X_m)}_ {m\in{\Zset}}$ is a  $\psi$-mixing process if for all integers $\ell\geq{0}$, the following
holds
\[\sup_{n\in{\Nset},B\in{{\mathcal{F}}_{\{0,.,n\}}},C\in{{\mathcal{F}}_{\{n\geq{0}\}}}}
\frac{|\p(B\cap{T^{-(n+\ell+1)}(C)})-\p(B)\p(C)|}
{\p(B)\p(C)}=\psi(\ell),\]
where the supremum is taken over the sets $B$ and $C$, such that $\p(B)\p(C)>0$.
\end{defn}

For a word $A$ of $\Omega$, that is to say a measurable subset of $\Omega$, we say that
$A\in{{\mathcal{C}}_n}$ if and only if 
\[A=\{X_0=a_0, \dots, X_{n-1}=a_{n-1}\}, \]
with $a_i \in \A, i=1, \dots, n$. Then, the integer $n$ is the length of word $A$.
For  $A\in{{\mathcal{C}}_n}$, we define the hitting time $\tau_{A}:\Omega
\rightarrow{\Nset\cup{\{\infty\}}}$, as the random variable defined
on the probability space ($\Omega$,$\mathcal{F}$,$\p$):
\[\forall x\in\Omega, \quad \tau_A(x)=\inf\{k\geq{1}:T^k(x)\in{A}\}.\]
$\tau_A$ is the first time that the process hits a given measurable $A$.
We also use the classical probabilistic shorthand notations. We write $\{\tau_{A}=m\}$ 
instead of $\{x\in{\Omega}:\tau_{A}(x)=m\}$, $T^{-k}(A)$ instead of
 $\{x\in{\Omega}:T^{k}(x)\in{A}\}$ and $\{X_r^s=x_r^s\}$ instead of 
$\{X_r=x_r,...,X_s=x_s\}$. Also we write for two measurable subsets $A$ and $B$ of $\Omega$, the
conditional probability of $B$ given $A$ as $\p(B|A)=\p_{A}(B)=\p(B \cap A)/\p(A)$ and the
probability of the intersection of $A$ and $B$ by $\p(A\cap B)$ or $\p(A;B)$.
For $A=\{X_0^{n-1}=x_0^{n-1}\}$ and  $1\leq{w}\leq{n}$, we write
$A^{(w)}=\{X_{n-w}^{n-1}=x_{n-w}^{n-1}\}$ for the event consisting of the \emph{last} $w$ symbols
of~$A$. We also write $a \vee b$ for the supremum of two real numbers $a$ and $b$. 
We define the periodicity $p_A$ of $A\in{{\mathcal{C}}_n}$ as follows:
\[p_A=\inf{\{k\in \Nset^*|A\cap T^{-k}(A)\neq {\emptyset}\}}.\]
$p_A$ is called the principal period of word $A$.
Then, we denote by $\mathcal{R}_p=\mathcal{R}_p(n)$ the set of words $A\in{\mathcal{C}_n}$ with 
periodicity $p$ and we also define $\mathcal{B}_n$ as the set of words  $A\in{\mathcal{C}_n}$
with periodicity less than~$\left[n/2\right]$, where $[.]$ defines the integer part of a real number:
\[\mathcal{R}_p=\{A\in{{\mathcal{C}}_n}|p_A=p\},   
\mathcal{B}_n=\bigcup_{p=1}^{\left[\frac{n}{2}\right]}\mathcal{R}_p.\]
$\mathcal{B}_n$ is the set of words which are self-overlapping before half their length (see
Example~\ref{mot}).
We define $\mathcal{R}(A)$ the set of return times of $A$ which are not a multiple of its
periodicity $p_A$:
\[\mathcal{R}(A)=\left\{ k \in \{[n/p_A]p_A+1, \dots, n-1\}| A\cap T^{-k}(A) 
\neq {\emptyset}\right\}.\]
Let us denote $r_A=\#\mathcal{R}(A)$, the cardinality of the set $\mathcal{R}(A)$. Define also
$n_A=\min{\mathcal{R}(A)}$ if $\mathcal{R}(A) \neq 
{\emptyset}$ and $n_A=n$ otherwise. $\mathcal{R}(A)$ is called the set of secondary periods of
$A$ and $n_A$ is the smallest secondary period of $A$.
Finally, we introduce the following notation. For an integer $s\in\{0, \dots, t-1\}$, let
$N_s^t=\sum_{i=s}^t \one\{T^{-i}(A)\}$. 
The random variable $N_s^t$ counts the number of occurrences of $A$ between $s$ and $t$ (we omit the
dependence on~$A$). For the
sake of simplicity, we also put $N^t = N_0^t$.
\begin{exmp}\label{mot}
Consider the word $A=aaataaataaa$. Since $p_A=4$, we have $A \in \mathcal{B}_n$ where $n=11$.
See the following figure to note that $\mathcal{R}(A)=\{9;10\}$, $r_A=2$ and $n_A=9$.
{\scriptsize
\begin{center}
$
\begin{array}{ccccccccccccccccccccc}
0&1&2&3&\mathbf{4}&5&6&7&\mathbf{8}&\mathbf{9}&\mathbf{10}&& & & & & & & & & \\
a&a&a&t&a&a&a&t&a&a&a& & & & & & & & & & \\
 & & & &a&a&a&t&a&a&a&t&a&a&a& & & & & & \\
 & & & & & & & &a&a&a&t&a&a&a&t&a&a&a& & \\
 & & & & & & & & &a&a&a&t&a&a&a&t&a&a&a& \\
 & & & & & & & & & &a&a&a&t&a&a&a&t&a&a&a\\
\end{array}
$
\end{center}
}
\end{exmp}

\subsection{The mixing method} \label{mixmet}

We present a theorem that gives an error bound for the Poisson approximation.
Compared to the Chen-Stein method, it has the advantage to present non uniform bounds 
that strongly control the  decay of the tail distribution of $N^t$.
% PSI MIXING
\begin{thm}[$\psi$-mixing approximation]\label{bigt2}
Let ${(X_m)}_{m\in{\Zset}}$ be a $\psi$-mixing process. There
 exists a constant $C_{\psi}=254$, such that for all $A\in{{\mathcal{C}}_n} \setminus{\mathcal{B}_n} $
  and
all non negative integers $k$ and $t$, the following inequality holds: 
\begin{eqnarray}
\left|\p(N^t=k)-\frac{{e}^{-t\p(A)}{(t\p(A))}^k}{k!}\right|
& \leq & C_{\psi} e_{\psi}(A) {e}^{-(t-(3k+1)n)\p(A)} g_{\psi}(A,k) \nonumber
\end{eqnarray}
\begin{eqnarray}
\mbox{ where } g_{\psi}(A,k) & = &
{
\left\{
\begin{array}{ll}
\displaystyle\frac{{(2\lambda)}^{k-1}}{(k-1)!} & 
\quad k\notin{\{\frac{\lambda}{e_{\psi}(A)},...,\frac{2t}{n}\}} \\
\displaystyle\frac{{(2\lambda)}^{k-1}}
{\left(\frac{\lambda}{e_{\psi}(A)}\right)!
{\left(\frac{1}{e_{\psi}(A)}\right)}^{k-\frac{\lambda}{e_{\psi}(A)}-1}}
& 
\quad k\in{\{\frac{\lambda}{e_{\psi}(A)},...,\frac{2t}{n}\}} 
\end{array}
\right.}, \nonumber
\end{eqnarray}
\[e_{\psi}(A)=\inf_{1\leq{w}\leq{n_A}}\left[(r_A+n)\p\left(A^{(w)}\right)\left(1+
\psi\left(n_A-w\right)\right)\right],\]
\[ \mbox{ and } \lambda={t\p(A)(1+\psi(n))}.\]
\end{thm}

This result is at the core of our study. It shows an upper bound for the difference between the
distribution of the number of occurrences of word $A$ in a sequence of length $t$ and the Poisson
distribution of parameter $t\p(A)$. 
Proof is postponed in Section \ref{Fin}.

\section{Calculation of the constants}\label{Prel}

Our goal is to compute a bound as small as possible to control the error between the
Poisson distribution and the distribution of the number of occurrences of a word.
Thus, we determine the global constant $C_{\psi}$ appearing in Theorem \ref{bigt2}  by means of
intermediary bounds appearing in the proof.  
General bounds are interesting asymptotically in $n$,
but for biological applications, $n$ is approximately between $10$ or $20$, which is too small.
Then along the proof, we will indicate the intermediary bounds that we compute.
Before establishing the proof of that Theorem \ref{bigt2}, we point out
here, for easy references, some results of \citet{abadi3}, 
and some other useful results. In \citet{abadi3}, these results are given only in the $\phi$-mixing
context. Moreover exact values of the constants are not given,
while these are necessary for practical use of these methods. We provide the values of all the
constants appearing in the proofs of these results.

%PROP 11
\begin{prop}[Proposition 11 in \citet{abadi3}]\label{p11}
Let ${(X_m)}_{m\in{\Zset}}$ be a $\psi$-mixing process.
There exist two finite constants $C_{a}>0$ and $C_{b}>0$, such that for any $n$, any word
 $A \in{{\mathcal{C}}_n} $,
and any $c\in \left[4n,\frac{1}{2\p\left(A\right)}\right]$ satisfying
\[\psi\left(c/4\right)\leq \p \left(\{\tau_A \leq c/4\} \cap \{\tau_A \circ T^{c/4}>c/2\} \right),\]
there exists $\Delta$, with $n<\Delta\leq c/4$, such that for all positive integers $k$,
the following inequalities hold:
\begin{eqnarray}
\left|\p\left(\tau_A>kc\right)-\p\left(\tau_A>c-2\Delta\right)^k\right| 
& \leq & C_{a} \varepsilon\left(A\right) k\p\left(\tau_A>c-2\Delta\right)^k, \label{in6} \\
\left|\p\left(\tau_A>kc\right)-\p\left(\tau_A>c\right)^k\right|
& \leq & C_{b} \varepsilon\left(A\right)k\p\left(\tau_A>c-2\Delta\right)^k, \label{in7} 
\end{eqnarray}
\[\mbox{with }\varepsilon(A)=\inf_{n\leq{{\ell}}\leq{\frac{1}{\p(A)}}}[{\ell}\p(A)+\psi({\ell})].\]
\end{prop}

%End PROP 11
Both inequalities provide an approximation of the hitting time distribution by a geometric
distribution at any point $t$ of the form 
$t=kc$. The difference between these distributions is that in \ref{in6}, the geometric term inside
the modulus is 
the same as in the upper bound, while in \ref{in7}, the geometric term inside the modulus is larger
than the one in the upper bound. That is, the second bound gives a larger error. We will use both in
the proof of Theorem \ref{ht}.

\begin{prop}
We have $C_a=24$ and $C_b=25$.
\end{prop}

\begin{pf}
For the details of the proof of Proposition \ref{p11}, we refer to Proposition $11$
in~\citet{abadi3}.
For any $c\in \left[4n,\frac{1}{2\p\left(A\right)}\right]$ and $\Delta \in \left[ n, c/4 \right]$,
  we denote 
$\n_j^i=\left\{\tau_A\circ T^{ic+j\Delta}>c-j\Delta\right\}$ and   
$\n=\left\{\tau_A>c-2\Delta\right\}$ for the sake of simplicity.
\citet{abadi3} obtains the following bound: 
\[\forall k\geq 2,
\left|\p\left(\tau_A>kc\right)-\p\left(\n\right)^k\right| \leq (a) + (b) + (c), \mbox{ with }\]
$(a) =  \displaystyle\sum_{j=0}^{k-2}\p\left(\n\right)^j
\left|\p\left(\tau_A>\left(k-j\right)c\right)-
         \p\left(\tau_A>\left(k-j-1\right)c ;  \n_2^{k-j-1}\right)\right|$,\\
$(b) =\displaystyle\sum_{j=0}^{k-2}\p\left(\n\right)^j \left|\p\left(\tau_A>\left(k-j-1\right)c ;  
\n_2^{k-j-1}\right)-\p\left(\tau_A>\left(k-j-1\right)c\right)\p\left(\n_2^0\right)\right|$,\\
$(c) = \p\left(\n\right)^{\left(k-1\right)}\left|\p\left(\tau_A>c\right)-\p\left(\n\right)\right|$.\\
First, for any measurable 
$B \in \mathcal{F}_{\{(\ell+1)c,(\ell+2)c+n-1\}}$, we have $\p\left(B\right)+\psi\left(\Delta\right)
\leq 3 \psi\left(\Delta\right)  \leq \frac{3}{2}\varepsilon\left(A\right)$.
We can also remark that $\p\left(\n\right) \geq 1/2$. Then, by iteration of the mixing property,
 we have the following inequality for all $\ell \in \Nset$:
\[\p\left(\bigcap_{i=0}^{\ell} \n_1^i ;  B\right)  \leq  6 \p\left(\n\right)^{\ell+1}
 \varepsilon\left(A\right).\]
We apply this bound in the inequalities (14) and (15) of \citet{abadi3} to get\\
$(a)\leq \displaystyle\sum_{j=0}^{k-2}\p\left(\n\right)^j\left(6\p\left(\n\right)^{k-j-2+1}
 \varepsilon\left(A\right)\right)
=6(k-1)\varepsilon\left(A\right)\p\left(\n\right)^{\left(k-1\right)}$,\\
$(b) \leq \displaystyle\sum_{j=0}^{k-2}\p\left(\n\right)^j\left(6\p\left(\n\right)^{k-j-2+1}
 \varepsilon\left(A\right)\right)
=6(k-1)\varepsilon\left(A\right)\p\left(\n\right)^{\left(k-1\right)}$.\\
We also have
$(c)\leq \p\left(\n\right)^{k-1}\p\left(\n ;  \tau_A \circ T^{c-2\Delta}
\leq 2\Delta \right)
 \leq \varepsilon\left(A\right) \p\left(\n\right)^{k-1}$.\\
We obtain (\ref{in6}):
$\left|\p\left(\tau_A>kc\right)-\p\left(\n\right)^k\right|
\leq 24k\varepsilon\left(A\right)\p\left(\n\right)^{k}$.\\
We deduce (\ref{in7}): 
$\left|\p\left(\tau_A>kc\right)-\p\left(\tau_A>c\right)^k\right|
\leq 25k\varepsilon\left(A\right)\p\left(\n\right)^{k}$.\\
Then, $C_{a}=24$ and $C_{b}=25$.
\end{pf}

%THM Hitting
\begin{thm}[Theorem 1 in \citet{abadi3}]\label{ht}
Let ${(X_m)}_{m\in{\Zset}}$ be a  $\psi$-mixing process.
Then, there exist constants $C_{h}>0$ and $0<\Xi_1<1\leq \Xi_2<\infty$, such that for all  $n\in{\Nset}$
and any
$A\in{{\mathcal{C}}_n}$, there exists $\xi_A\in{[\Xi_1,\Xi_2]}$, for which the following inequality 
holds for all $t>0$:
\[\left|\p\left(\tau_{A}>\frac{t}{\xi_A}\right)-
{e}^{-t\p(A)}\right|\leq{C_{h}\varepsilon(A)f_1(A,t)},\]
\[\mbox{with }\varepsilon(A)=\inf_{n\leq{{\ell}}\leq{\frac{1}{\p(A)}}}[{\ell}\p(A)+\psi({\ell})]
\mbox{ and }
f_1(A,t)=(t\p(A)\vee1){e}^{-t\p(A)}.\]
\end{thm}
We prove an upper bound for the distance between the rescaled hitting time and the exponential law
of expectation equal to one. 
The factor $\varepsilon(A)$ in the upper bound shows that the rate of convergence to the exponential law
is given by a trade off between the length of this time and the velocity of loosing memory of the process.
%End THM Hitting

\begin{prop}
We have $C_h=105$.
\end{prop}

\begin{pf}
We fix $c=\frac{1}{2\p(A)}$ and $\Delta$ given by Proposition \ref{p11}. We define
\[\xi_A=\frac{-\log{\p(\tau_A>c-2\Delta)}}{c\p(A)}.\]
There are three steps in the proof of the theorem. First, we consider $t$ of the form $t=kc$ with $k$
 a positive integer. Secondly, we prove the theorem for any $t$ of the form $t=(k+p/q)c$ with $k,p$ 
positive integers
and $1\leq p\leq q$ with $q=\frac{1}{2\varepsilon(A)}$. We also put $r=(p/q)c$. Finally, 
 we consider the remaining cases.
Here, for the sake of simplicity , we do not detail the two first steps (for that, see
\citet{abadi3}), but only the last one. 
Let $t$ be any positive real number. We write $t=kc+r$, with $k$ a positive integer and $r$ such that
 $0\leq r <c$.
 We can choose a $\bar{t}$ such that $\bar{t}<t$ and $\bar{t}=(k+p/q)c$ with $p$, $q$ as before.
 \citet{abadi3} obtains the following bound:
\begin{eqnarray}
\left|\p\left(\tau_A>t\right)-{e}^{-\xi_A\p\left(A\right)t}\right| & \leq & 
\left|\p\left(\tau_A>t\right)-\p\left(\tau_A>\bar{t}\right)\right| 
 +  \left|\p\left(\tau_A>\bar{t}\right)-{e}^{-\xi_A\p\left(A\right)\bar{t}}\right| \nonumber \\
 & + & \left|{e}^{-\xi_A\p\left(A\right)\bar{t}}-{e}^{-\xi_A\p\left(A\right)t}\right|. \nonumber 
\end{eqnarray}
The first term in the triangular inequality is bounded in the following way:
\begin{eqnarray}
\left|\p\left(\tau_A>t\right)-\p\left(\tau_A>\bar{t}\right)\right|  & = &
\p\left( \tau_A >\bar{t} ;  \tau_A\circ T^{\bar{t}} \leq t-\bar{t}\right) \nonumber \\
 & \leq & \p\left( \tau_A > kc ;  \tau_A\circ T^{\bar{t}} \leq \Delta \right) \nonumber \\
 & \leq & {\p\left(\n\right)}^{k-2}\left(\Delta\p(A)+\psi(\Delta))\right) \nonumber \\
 & \leq & 4 {\p\left(\n\right)}^{k}\varepsilon(A) \nonumber \\
 & \leq & 4\varepsilon(A){e}^{-\xi_A\p(A)t}. \nonumber 
\end{eqnarray}

The second term is bounded like in the two first steps of the proof in \citet{abadi3}.
We apply inequalities (\ref{in6}) and (\ref{in7}) to obtain
\[\left|\p\left(\tau_A>\bar{t}\right)-{e}^{-\xi_A\p\left(A\right)\bar{t}}\right|
\leq  (3+C_{a} t\p(A) + C_{a} + 2C_{b})\varepsilon(A){e}^{-\xi_A\p(A)t}. \]
Finally, the third term is bounded using the Mean Value Theorem (see for example \citet{MVT})
\[\left|{e}^{-\xi_A\p\left(A\right)\bar{t}}-{e}^{-\xi_A\p\left(A\right)t}\right| 
\leq  \xi_A\p(A)\left(r-\frac{p}{q}c\right){e}^{-\xi_A\p(A)\bar{t}} 
\leq   \varepsilon(A){e}^{-\xi_A\p(A)t}.\]
Thus we have 
$\left|\p\left(\tau_A>t\right)-{e}^{-\xi_A\p\left(A\right)t}\right| \leq 
105 \varepsilon(A)f_1(A,\xi_A t)$
and the theorem follows by the change of variables $\tilde{t}=\xi_At$.
Then $C_{h}=105$.
\end{pf}

%Begin Lemma p(bc)
\begin{lem}\label{pbc}
$\left( X_m \right)_{m \in \Zset}$ be a $\psi$-mixing  process.
Suppose that $B \subseteq A \in \mathcal{F}_{\{0,\dots, b\}}$, $C \in \mathcal{F}_{\{b+g,\dots,
  \infty\}}$ with 
$b, g \in \Nset$. The following inequality holds:
\[\p_A(B\cap C) \leq \p_A(B)\p(C)(1+\psi(g)). \]
\end{lem}
%End Lemma p(bc)

\begin{pf}
Since $B \subseteq A$, obviously $\p(A\cap B \cap C)=\p(B \cap C)$. 
By the $\psi$-mixing 
property $\p(B\cap C) \leq \p(B)(\p(C)+\psi(g))$.
We divide the above inequality by $\p(A)$ and the lemma follows. 
\end{pf}

For all the following propositions and lemmas, we recall that 
$$e_{\psi}(A)=\inf_{1\leq{w}\leq{n_A}}\left[(r_A+n)\p\left(A^{(w)}\right)\left(1+
\psi\left(n_A-w\right)\right)\right].$$
%Begin Proposition 11 Stat
\begin{prop} \label{:lambda}
Let   $\left( X_m \right)_{m \in \Zset}$ be a $\psi$-mixing  process. Let $A \in {\mathcal R}_p(n)$. Then
the following holds:
\begin{itemize}
\item[(a)] For all $M, M' \geq g\geq n$, 
\begin{eqnarray}
& & \left\vert \p_A\left(\tau_A>M+M'\right)- \p_A\left(\tau_A> M \right) \p\left( \tau_A> M' \right)
\right\vert \nonumber \\
&\leq& \p_A\left( \tau_A> M-g \right) 2g\p(A) \left[ 1+\psi(g) \right],\nonumber 
\end{eqnarray}
and similarly
\begin{eqnarray}
& & \left\vert\p_A\left(\tau_A>M + M'\right) -
  \p_A\left(\tau_A>M\right)\p\left(\tau_A>M'-g\right)\right\vert \nonumber \\
&\leq& \p_A\left(\tau_A>M-g\right)\left[ g\p(A)+2\psi(g) \right].\nonumber 
\end{eqnarray}
\item[(b)] For all $t \geq p \in \Nset$, with $\zeta_A=\p_A(\tau_A>p_A)$,
\[
\left\vert \p_A\left( \tau_A> t \right) - \zeta_A \p\left( \tau_A> t \right) \right\vert \leq 2 e_{\psi}(A).
\]
\end{itemize}
\end{prop}
%End Proposition 11 Stat

The above proposition establishes a relation between hitting and 
return times with an error bound uniform with respect to $t$. In particular, $(b)$ says that these times coincide if and only if  
$\zeta_A= 1$, namely, the string $A$ is non-self-overlapping. 

\begin{pf} 
In order to simplify notation, for $t\in \Zset$, $\tau_A^{[t]}$ stands for $\tau_A\circ T^t$.
We introduce a gap of length $g$
after coordinate $M$ to construct the following triangular
inequality
\begin{eqnarray}
&&
\left\vert \p_A\left( \tau_A> M+M' \right) - \p_A\left( \tau_A> M \right) \p\left( \tau_A> M' \right)
\right\vert \nonumber\\ 
&\leq&
\left\vert\p_A\left( \tau_A> M+M' \right) -
   \p_A\left( \tau_A> M  ;  \tau_A^{[M+g]}>M'-g \right)\right\vert \label{mix1} \\
&+&
\left\vert\p_A\left(\tau_A>M ; \tau_A^{[M+g]}>M'-g\right)-
   \p_A\left(\tau_A>M\right)\p\left(\tau_A>M'-g\right)\right\vert \label{mix2} \\
&+& 
    \p_A\left( \tau_A> M \right)
\left\vert \p\left( \tau_A> M'-g \right) -  \p\left(\tau_A> M' \right) \right\vert.\label{mix3} 
\end{eqnarray}
Term (\ref{mix1}) is bounded with Lemma~\ref{pbc} by
\[
    \p_A\left( \tau_A >  M  ;  \tau_A^{[M]} \leq g\right)
\leq \p_A\left( \tau_A> M -g \right) g\p(A)\left[ 1+\psi(g) \right].
\]
Term (\ref{mix2}) is bounded using the $\psi$-mixing property by
$ \p_A\left( \tau_A> M \right) \psi(g)  $.
The modulus in (\ref{mix3}) is bounded using stationarity by
$\p\left( \tau_A \leq g \right) \leq g \p(A) .$
This ends the proof of both inequalities of item $(a)$.\\
Item $(b)$ for $t\geq 2n$ is proven similarly to item $(a)$ with
$t=M+M'$, $M=p$, and $g=w$ with $1\leq w \leq n_A$. Consider now
$p\leq  t < 2n$.
\[
   \zeta_A -  \p_A\left( \tau_A > t \right)
= \p_A\left( p < \tau_A \leq t \right)
 =  \p_A\left( \tau_A \in \mathcal{R}(A) \cup \left(  n \leq \tau_A \leq t \right) \right) \leq
 e_{\psi}(A). \nonumber 
\]
First and second equalities follow by definition of $\tau_A$ and $\mathcal{R}(A)$. The inequality
follows by Lemma \ref{pbc}.
\end{pf}

Let $\zeta_A = \p_A(\tau_A>p_A)$ and $h=1/(2\p(A))-2\Delta$, 
then $\xi_A = -2\log{\p(\tau_A>h)}$.
%Lemma 12
\begin{lem}\label{l12}
Let  ${(X_m)}_{m\in{\Zset}}$ be a $\psi$-mixing  process. Then the
following inequality holds:
\[|\xi_A-\zeta_A| \leq 11 e_{\psi}(A).\]
\end{lem}
%End Lemma 12
Hence, we have
\[\zeta_A-11e_{\psi}(A)\leq \xi_A \leq \zeta_A+11e_{\psi}(A).\]

\begin{pf}
\begin{eqnarray}
\p\left( \tau_A  > h \right)
& = &\prod_{i=1}^{h} \p\left( \tau_A>i | \tau_A  > i-1 \right) 
=\prod_{i=1}^{h} (1-\p\left( T^{-i}(A) | \tau_A  > i-1 \right)) \nonumber \\
& = &\prod_{i=1}^{h} \left( 1- \rho_i \p(A) \right),\nonumber 
\end{eqnarray}
where
$
\rho_i
\stackrel{def}{=}
\displaystyle\frac{\p_A\left( \tau_A  > i-1 \right)}{\p\left( \tau_A>i-1 \right)}.
$
Therefore
\begin{eqnarray}
& & \left\vert \xi_A + 2\sum_{i=1}^{p_A}\log(1-\rho_i\p(A)) 
        - 2\sum_{i=p_A+1}^{h}\zeta_A\p(A) \right\vert \nonumber \\
&\leq& 
2\sum_{i=p_A+1}^{h} 
\left\vert -\log(1-\rho_i\p(A)) -  \zeta_A \p(A) \right\vert.\nonumber 
\end{eqnarray}
The above modulus is bounded by
\[
\left\vert -\log(1-\rho_i\p(A)) -  \rho_i \p(A) \right\vert
+
\left\vert \rho_i -  \zeta_A \right\vert \p(A).
\]
Now note that $|y - (1 - e^{-y}) | \leq (1- e^{-y})^2$ for $y>0$ small enough.
Apply it with $y=-\log(1-\rho_i\p(A))$ to bound
the most left term of the above expression by $(\rho_i\p(A))^2$.
Further by Proposition \ref{:lambda} $(b)$ and the
fact that $\p\left(\tau_A>h\right)\geq 1/2$ we have
\[
\left\vert \rho_i - \zeta_A \right\vert \leq \frac{2 e_1(A)}{\p\left( \tau_A > h \right)}
\leq 4 e_{\psi}(A).
\]
for all $i=p_A+1,\dots,h$. 
Yet as before
\[
- \sum_{i=1}^{p_A}\log(1-\rho_i\p(A)) 
\leq p_A \left( \rho_i\p(A) + (\rho_i\p(A))^2 \right) \leq  e_{\psi}(A).
\]
Finally, by definition of $h$
\[
\left\vert 2\sum_{i=p_A+1}^{h} \zeta_A \p(A)  - \zeta_A \right\vert
\leq 4\Delta\p(A) + 2 p_A\p(A) \leq 6 e_{\psi}(A).
\]
This ends the proof of the lemma.
\end{pf}

%PROP
\begin{prop}\label{mprop}
Let ${(X_m)}_{m\in{\Zset}}$ be a  $\psi$-mixing process. Then the following
inequality holds:
\[ |\p(\tau_A>t)-{e}^{-t\p(A)}| \leq C_{p}e_{\psi}(A)(t\p(A)\vee1){e}^{-(\zeta_A-11e_{\psi}(A))t\p(A)}. \]
\end{prop}
%End PROP

\begin{pf}
We bound the first term with Theorem \ref{ht} and the second with Lemma~\ref{l12}~:
\[
\begin{array}{ccl}
|\p(\tau_A>t)-{e}^{-t\p(A)}| & \leq &
|\p(\tau_A>t)-{e}^{-\xi_A t\p(A)}|+|{e}^{-\xi_A t\p(A)}-{e}^{-t\p(A)}| \\
|\p(\tau_A>t)-{e}^{-\xi_A t\p(A)}|  & \leq & 
C_{h} \varepsilon(A){e}^{-\xi_A t\p(A)} \leq  
C_{h} e_{\psi}(A){e}^{-(\zeta_A-11e_{\psi}(A))t\p(A)} \\
|{e}^{-\xi_A t\p(A)}-{e}^{-t\p(A)}| & \leq & 
t\p(A)|\xi_A-1|{e}^{-\min{\{1,\xi_A\}t\p(A)}} \\
 & \leq & 11t\p(A) e_{\psi}(A){e}^{-(\zeta_A-11e_{\psi}(A))t\p(A)}. \\
\end{array}
\]
This ends the proof of the proposition with $C_{p} = C_{h}+11$.
\end{pf}

%DEF
\begin{defn}
Given $A\in{{\mathcal{C}}_n}$, we define for $j\in{\Nset}$, the $j$-th occurrence time
of  $A$ as the random variable $\tau_A^{(j)}:\Omega\rightarrow{\Nset\cup\{\infty\}}$,
defined on the probability space $(\Omega, \mathcal{F}, \p)$ as follows: for any 
$x\in{\Omega}$, $\tau_A^{(1)}(x)=\tau_A(x)$ and for $j\geq{2}$,
\[\tau_A^{(j)}(x)=\inf{\{k>\tau_A^{(j-1)}(\omega):T^k(x)\in{A}\}}.\]
\end{defn}
%End DEF

%PROP 18
\begin{prop}\label{18}
Let ${(X_m)}_{m\in{\Zset}}$ be a $\psi$-mixing process. Then, for all $A\notin{\mathcal{B}_n}$, all $k \in \Nset$,
and all $0\leq{t_1}<t_2<...<t_k\leq{t}$ for which
$\displaystyle\min_{2\leq{j}\leq{k}}\{t_j-t_{j-1}\}>2n$, there
 exists a positive constant $C_{1}$ independent of $A$, $n$, $t$ and $k$ such that
\begin{eqnarray}
& & \left|\p\left(\displaystyle \bigcap_{j=1}^k\left(\tau_A^{(j)}=t_j\right) ; 
 \tau_A^{(k+1)}>t\right)-
{\p(A)}^k\displaystyle \prod_{j=1}^{k+1}\mathcal{P}_j\right| \nonumber \\
& \leq& {C_{1}k{(\p(A)(1+\psi(n)))}^k e_{\psi}(A)}{{e}^{-(t-(3k+1)n)\p(A)}} \nonumber
\end{eqnarray}
where $\mathcal{P}_j=\p(\tau_A>(t_j-t_{j-1})-2n)$.
\end{prop}
%End PROP 18

\begin{pf}
We will show this proposition by induction on $k$.
We put  $\Delta_j=t_j-t_{j-1}$ for $j=2,...,k$, $\Delta_1=t_1$ and $\Delta_{k+1}=t-t_k$. 
Firstly, we note that by stationarity \[\p(\tau_A=t)=\p(A   ;    \tau_A>t-1).\]
For $k = 1$, by a triangular inequality we obtain
\begin{eqnarray}
&      & \left|\p\left(\tau_A=t_1   ;    \tau_A^{(2)}>t\right)-
\p(A)\prod_{j=1}^{2}\mathcal{P}_j\right| \nonumber \\
& \leq & \left|\p\left(\tau_A=t_1   ;    \tau_A^{(2)}>t\right)-
\p\left(\tau_A=t_1   ;   N_{t_1+2n}^t=0\right)\right| \label{18a} \\
& +    &\left|\p\left(\tau_A=t_1   ;   N_{t_1+2n}^t=0\right)-
\p\left(\tau_A=t_1\right)\mathcal{P}_2\right| \label{18b} \\
& +    &\left|\p(A   ;    \tau>t_1-1)-
\p\left(A   ;   N_{2n}^{t_1-1}=0\right)\right|\mathcal{P}_2 \label{18c} \\
& +    &\left|\p\left(A   ;   N_{2n}^{t_1-1}=0\right)\mathcal{P}_2-
\p(A)\prod_{j=1}^{2}\mathcal{P}_j\right|. \label{18d}
\end{eqnarray}
Term (\ref{18a}) is equal to
$\p\left(\tau_A=t_1 ;  \bigcup_{i=t_1+1}^{t_1+2n}T^{-i}(A) ;  N_{t_1+2n}^{t}=0 \right)$ and then

$$(\ref{18a})=\p\left(A ;  \bigcup_{i\in\mathcal{R}(A) \cup i=1}^{2n}T^{-i}(A) ;  N_{2n}^{t}=0 \right).$$
Since $A \notin \mathcal{B}_n$, for $1\leq i <p_A$, the above probability is zero. Thus, using
mixing property 
\begin{eqnarray}
(\ref{18a}) & \leq & \p\left(A ;  \bigcup_{i\in\mathcal{R}(A) \cup i=p_A}^{2n}T^{-i}(A) ; 
  N_{2n}^{t}=0 \right) \nonumber\\
 & \leq & 2\p(A)\p(A)(r_A+n)(1+\psi(n)) \p\left(N_{2n}^{t}=0\right) \nonumber \\
 & \leq & 2\p(A) e_{\psi}(A) e^{-(t-(3k+1)n)\p(A)}. \nonumber 
\end{eqnarray} 
Term (\ref{18b}) is bounded using $\psi$-mixing property
\begin{eqnarray}
(\ref{18b}) & \leq & \psi(n)(1+\psi(n))\p(A)\mathcal{P}_1\mathcal{P}_2 \nonumber \\
 & \leq & \psi(n)\p(A)e_{\psi}(A)e^{-(t-(3k+1)n)\p(A)}. \nonumber 
\end{eqnarray}
Analogous computations are used to bound terms (\ref{18c}) and (\ref{18d}).

Now, let us suppose that the proposition holds for $k-1$ and let us prove it for $k$.
We put $\mathcal{S}_i =\{\tau_A^{(i)}=t_i\}$.
We use a triangular inequality again to bound the term in the left hand side
of the inequality of the proposition by a sum of five terms:
\newline
$
\left|\p\left(\displaystyle \bigcap_{j=1}^k\left(\tau_A^{(j)}=t_j\right)   ;   
      \tau_A^{(k+1)}>t\right)-
    {\p(A)}^k\displaystyle \prod_{j=1}^{k+1}\mathcal{P}_j\right| 
\leq  I+II+III+IV+V.
$
\newline
$
\begin{array}{lll}
  I & = &  \left|\p\left(\displaystyle\bigcap_{j=1}^k\mathcal{S}_j   ;   
      \tau_A^{(k+1)}>t\right) 
  \right.
   -  \left. \p\left(\displaystyle\bigcap_{j=1}^{k-1}\mathcal{S}_j   ;   
      N_{t_{k-1}+1}^{t_k-2n}=0   ;   T^{-t_k}(A)   ;   N_{t_k+1}^t=0\right)\right|\nonumber \\
  & = & \p\left(\displaystyle\bigcap_{j=1}^{k-1}\mathcal{S}_j ; 
  N_{t_{k-1}+1}^{t_k-2n}=0 ;  \displaystyle\bigcup_{i=t_k-2n+1}^{t_k-1}T^{-i}(A) ;  T^{-t_k}(A) ; 
  N_{t_{k}+1}^{t}=0 \right) \nonumber \\ 
  & \leq &{(\p(A)(1+\psi(n)))}^k(1-\psi(n))\left(np_A+ {(r_A+n)\p(A^{(w)})}\right)
  e^{-(t-(3k+1)n)\p(A)}, \nonumber 
\end{array}
$
\newline
$
\begin{array}{lll}
  II& = & \left|\p\left(\displaystyle\bigcap_{j=1}^{k-1}\mathcal{S}_j   ;   
      N_{t_{k-1}+1}^{t_k-2n}=0   ;   T^{-t_k}(A)   ;   N_{t_k+1}^t=0\right)\right.\nonumber \\
  & - & \left. \p\left(\displaystyle\bigcap_{j=1}^{k-1}\mathcal{S}_j   ;   
      N_{t_{k-1}+1}^{t_k-2n}=0\right)\p\left(A   ;   N_1^{t-t_k}=0\right)\right|\nonumber\\
  &\leq& \p\left(\displaystyle\bigcap_{j=1}^{k-1} ;  N_{t_{k-1}+1}^{t_k-2n}=0\right)
     \p\left(A ;  N_1^{t-t_k}=0\right)\psi(n) \\
  &\leq& {(\p(A)(1+\psi(n)))}^k\psi(n)e^{-(t-(3k+1)n)\p(A)},\nonumber 
\end{array}
$
\newline
$
\begin{array}{lll}
  III & =& \left|\p\left(\displaystyle\bigcap_{j=1}^{k-1}\mathcal{S}_j   ;   
      N_{t_{k-1}+1}^{t_k-2n}=0\right)\right. 
   -  \left. \p\left(\displaystyle\bigcap_{j=1}^{k-1}\mathcal{S}_j   ;   
      N_{t_{k-1}+1}^{t_k-1}=0\right) \right| \p\left(A    ;   N_1^{t-t_k}=0\right)  \nonumber\\
  &\leq& \p\left(\displaystyle\bigcap_{j=1}^{k-1}\mathcal{S}_j ;  
      N_{t_{k-1}+1}^{t_k-2n}=0 ;  \displaystyle\bigcup_{t_k-2n+1}^{t_k-1}T^{-i}(A)\right)\p(A)  \nonumber\\
  &\leq& 2\p(A){(\p(A)(1+\psi(n)))}^k e^{-(t-(3k+1)n)\p(A)}. \nonumber
\end{array}
$
\newline
We use the inductive hypothesis for the term $IV$ and the case with $k=1$ for the term~$V$. 
\newline
$
\begin{array}{lll}
  IV & =      &  \left|\p\left(\displaystyle\bigcap_{j=1}^{k-1}\mathcal{S}_j ; 
      N_{t_{k-1}+1}^{t_k-1}=0\right)-{\p(A)}^{k-1}\displaystyle\prod_{j=1}^k\mathcal{P}_j\right|
  \p\left(A ; N_1^{t-t_k}=0\right)\nonumber \\
  & \leq   & C_{1}(k-1){(\p(A)(1+\psi(n)))}^k e_{\psi}(A)e^{-(t-(3k+1)n)\p(A)}, \nonumber
\end{array}
$
\newline
$
\begin{array}{lll}
  V & = & {\p(A)}^{k-1}
  \displaystyle\prod_{j=1}^k\mathcal{P}_j\left |\p\left(A   ;   N_1^{t-t_k}=0\right)-
    \p(A)\mathcal{P}_{k+1} \right|  \nonumber \\
  &\leq& 2{(\p(A)(1+\psi(n)))}^k e_{\psi}(A)e^{-(t-(3k+1)n)\p(A)}. \nonumber
\end{array}
$
\newline
Finally, we obtain \[I+II+III+IV+V\leq (3+C_{1}(k-1)+2){(\p(A)+\psi(n))}^ke_{\psi}(A).\]
To conclude the proof, it is sufficient that $C_{1}k=3+C_{1}(k-1)+2$, therefore $C_{1}=5$.
This ends the proof of the proposition.
\end{pf}

\section{Proof of Theorem \ref{bigt2}}\label{Fin}

In this section, we prove the main result of our work (see Section \ref{mixmet}): an upper bound for the
difference between the 
exact distribution of the number of occurrences of word $A$ and the Poisson distribution of
parameter $t\p(A)$.
Throughout the proof, we will note in italic the terms computed by our software \texttt{PANOW} (see Section
\ref{PANOW}).

\begin{pf}
For $k=0$, the 
result comes from Proposition~\ref{mprop} ($\p(N^t=0)=\p(\tau_A>t)$).
\\
For $k>2t/n$, since $A\notin{\mathcal{B}_n}$, we have $\p(N^t=k)=0$. Hence,
\begin{eqnarray}
\left|\p(N^t=k)-\frac{{e}^{-t\p(A)}{(t\p(A))}^k}{k!}\right| 
 & = &\frac{{e}^{-t\p(A)}{(t\p(A))}^k}{k!} \nonumber \\
 & \leq &{\frac{(t\p(A))^{k-1}}{(k-1)!}\frac{t\p(A)}{k}} \nonumber\\ 
 & \leq & \frac{1}{2}\frac{(t\p(A))^{k-1}}{(k-1)!}e_{\psi}(A). \nonumber
\end{eqnarray}
Indeed, since $\frac{t}{k}<\frac{n}{2}$ then 
$\frac{t\p(A)}{k} < \frac{n\p(A)}{2} \leq {\frac{e_{\psi}(A)}{2}}$. 

Now, let us consider $1\leq{k}\leq{2t/n}$.
We consider a sequence which contains exactly 
$k$ occurrences of $A$. These occurrences can be isolated or can be in clumps. We define the
 following set:
\[\mathcal{T} =\mathcal{T}( t_1,t_2,...,t_k)=\left\{ \bigcap_{j=1}^k(\tau_A^{(j)}=t_j)   ;   
\tau_A^{(k+1)}>t\right\}.\]
We recall that we put $\mathcal{P}_j=\p(\tau_A>(t_j-t_{j-1})-2n)$,
 $\Delta_j=t_j-t_{j-1}$ for $j=2,...,k$, $\Delta_1=t_1$ and $\Delta_{k+1}=t-t_k$. 
Define $I(\mathcal{T})=\displaystyle \min_{2\leq{j}\leq{k}}\{\Delta_j\}$.
We say that the occurrences of $A$ are isolated if $I(\mathcal{T})\geq 2n$ and we say that
there exists at least one clump if $I(\mathcal{T})<2n$.
We also denote 
\[B_k=\left\{\mathcal{T}  |  I(\mathcal{T})< 2n \right\} \quad \mbox{ and } \quad
  G_k=\left\{\mathcal{T}  |  I(\mathcal{T})\geq 2n \right\}.\]
The set $\{N^t=k\}$ is the disjoint union between $B_k$ and $G_k$, then \[\p(N^t=k)=\p(B_k)+\p(G_k),
\]
\[
\left|\p(N^t=k)-\frac{{e}^{-t\p(A)}{(t\p(A))}^k}{k!}\right| \leq{\p(B_k)}+
\left|\p(G_k)-\frac{{e}^{-t\p(A)}{(t\p(A))}^k}{k!}\right|.
\]
We will prove an upper bound for the two quantities on the right hand side of the above inequality
to conclude the proof of the theorem.

\textbf{We prove an upper bound for} $\mathbf{\p(B_k)}$\textbf{.}  Define 
$C(\mathcal{T}) = \sum_{j=2}^{k} {\one}_{\left\{ \Delta_{j} > 2n \right\}} + 1$. 
$C(\mathcal{T})$ computes how many clusters there
are in a given $\mathcal{T}$.
Suppose that $\mathcal{T}$ is such that
$C(\mathcal{T}) = 1$ and fix the position $t_1$ of the first occurrence
of $A$. Further, each occurrence inside the cluster (with the
exception of the most left one which is fixed at $t_1$) can appear at
distance $d$ of the previous one, with $p_A \leq d \leq 2n$. Therefore, 
the $\psi$-mixing property leads to the bound
\begin{eqnarray}
\p\left(\bigcup_{t_2,\dots,t_k}
\mathcal{T}(t_1,t_2,\dots,t_k) \right) &\leq& \p\left( \bigcap_{j=1}^{k} 
      \bigcup_{
\begin{array}{c}        
\scriptstyle n/2 \leq t_{i+1}-t_i \leq 2n; \\
\scriptstyle i=2,\dots,k 
\end{array}
} 
T^{-t_j}(A)  \right)  \\
&\leq& \p(A) e_{\psi}(A)^{k-1}e_{\psi}(A)e^{-(t-(3k+1)n)\p(A)} . \nonumber
\end{eqnarray}
Suppose now that $\mathcal{T}$ is such that $C(\mathcal{T})=i$.
Assume also that
the most left occurrence of the $i$ clusters of $\mathcal{T}$
occurs at $t(1),\dots,t(i)$, with $1\leq t(1)<\dots<t(i)\leq t$ fixed.
By the same argument used above, we have the inequalities
\begin{eqnarray}
&  & \p\left( \bigcup_{ \{t_1,\dots,t_k\} \backslash \{t(1),\dots,t(i)\} }
\mathcal{T}(t_1,\dots,t_k) \right) \nonumber \\
&\leq& {\left(\p(A)(1+ \psi(n))\right)}^{i-1} e_{\psi}(A)^{k-i}e^{-(t-(3k+1)n)\p(A)}. \nonumber
\end{eqnarray}
To obtain an upper bound for $ \p\left(B_{k}\right)$ we must sum the above bound
over all $\mathcal{T}$ such that $C(\mathcal{T} ) = i$ with $i$
running from $1$ to $k-1$.
Fixed $C(\mathcal{T})=i$, the locations of the most left occurrences of $A$ of each one
of the $i$ clusters can be chosen in at most  $C_t^i$ many ways.
The cardinality of each one of the $i$ clusters can be arranged in
$C_{k-1}^{i-1}$ many ways. (This corresponds to breaking the
interval $(1/2,k+1/2)$ in $i$ intervals at points chosen  from
$\{1+1/2,\dots,k-1/2\}$.)
Collecting these informations, we have that $ \p\left(B_{k}\right)$ is bounded by
\[
\begin{array}{lll}
  &     &  {\displaystyle\sum_{i=1}^{k-1}C_t^i C_{k-1}^{i-1}
  {(\p(A)(1+\psi(n)))}^i e_{\psi}(A)^{k-i}} {e}^{-(t-(3k+1)n)\p(A)} \\
  &\leq &  {e}^{-(t-(3k+1)n)\p(A)} e_{\psi}(A)^{k}\displaystyle\max_{1\leq i \leq k-1}
 \displaystyle\frac{{\left(\lambda / e_{\psi}(A)\right)}^i}{i!}
  \displaystyle\sum_{i=1}^{k-1}C_{k-1}^{i-1} \\ 
  & \leq & {e}^{-(t-(3k+1)n)\p(A)} e_{\psi}(A){ \left\{
      \begin{array}{ll}
        \displaystyle\frac{{(2\lambda)}^{k-1}}{(k-1)!} &
        k<\frac{\lambda}{e_{\psi}(A)} \\ 
        \displaystyle\frac{{(2\lambda)}^{k-1}}{\left(\frac{\lambda}{e_{\psi}(A)}\right)!
          {\left(\frac{\lambda}{e_{\psi}(A)}\right)}^{k-1-\frac{\lambda}{e_{\psi}(A)}}}&
        k\geq \frac{\lambda}{e_{\psi}(A)} \\
      \end{array}
    \right.}.\\
\end{array}
\]
This ends the proof of the bound for $\p\left( B_{k}\right)$.
\\
{\itshape
We compute $\p(B_k) \leq {\displaystyle\sum_{i=1}^{k-1}C_t^i C_{k-1}^{i-1}
  {(\p(A)(1+\psi(n)))}^i e_{\psi}(A)^{k-i}} {e}^{-(t-(3k+1)n)\p(A)}$. 
}

\textbf{We prove an upper bound for} 
$\left|\p(G_k)-\displaystyle\frac{{e}^{-t\p(A)}{(t\p(A))}^k}{k!}\right|$\textbf{.} 
It is bounded by four terms by the triangular inequality 
\begin{eqnarray}
 & & \sum_{T\in{G_k}}\left|\p\left( \bigcap_{j=1}^k\left(\tau_A^{(j)}=t_j\right)   ;   
\tau_A^{(k+1)}>t\right)
-{\p(A)}^k\prod_{j=1}^{k+1}{\mathcal{P}}_j\right| \label{14}\\
 &+& \sum_{T\in{G_k}}{\p(A)}^k\left|\prod_{j=1}^{k+1}{\mathcal{P}}_j
-\prod_{j=1}^{k+1}{e}^{-\left(\Delta_j-2n\right)\p(A)}\right| \label{15}\\
 &+& \sum_{T\in{G_k}}{\p(A)}^k\left|{e}^{-\left(t-2(k+1)n\right)\p(A)}
-{e}^{-t\p(A)}\right| \label{16} \\
 &+& \left|\frac{\#G_k  k!}{t^k}\frac{{e}^{-t\p(A)}{(t\p(A))}^k}{k!}
-\frac{{e}^{-t\p(A)}{(t\p(A))}^k}{k!}\right|.\label{17}
\end{eqnarray}
We will bound these terms to obtain Theorem \ref{bigt2}.\newline
First, we bound the cardinal of $G_k$ 
\[
\#G_k
\leq
C_t^k
\leq
{\frac{t^k}{k!}}.
\]
Term (\ref{14}) is bounded with Proposition \ref{18} 
\[(\ref{14})\leq{C_{1}\frac{t^k}{(k-1)!}{(\p(A)(1+\psi(n)))}^k e_{\psi}(A)}{e}^{-(t-(3k+1)n)\p(A)}.\]
Term (\ref{15}) is bounded with Proposition \ref{mprop} 
\begin{eqnarray}
(\ref{15})& \leq & \frac{t^k}{k!}{\p(A)}^k\sum_{j=1}^{k+1}\prod_{i=1}^{j-1}{\mathcal{P}}_i
\left|{\mathcal{P}}_j-{e}^{-\left(\Delta_j-2n\right)\p(A)}\right|\prod_{i=j+1}^{k+1}
{e}^{-\left(\Delta_i-2n\right)\p(A)} 
\nonumber\\
 & \leq & \frac{t^k}{k!}{\p(A)}^k(k+1)   C_{p}e_{\psi}(A)  {e}^{-(\zeta_A-11e_{\psi}(A)) t\p(A)}\nonumber\\
 & \leq & 2  C_{p}\frac{{(t\p(A))}^k}{(k-1)!}e_{\psi}(A) {e}^{-(\zeta_A-11e_{\psi}(A)) t\p(A)}\nonumber
\end{eqnarray}
where $C_{p}$ is defined in Proposition \ref{mprop}. 
\\
{\itshape
We compute
\begin{eqnarray}
(\ref{15}) & \leq & \frac{{(t\p(A))}^k}{(k-1)!}\frac{k+1}{k} \nonumber\\
          &  & \left[(8+C_{a} t\p(A) + C_{a} + 2C_{b})\varepsilon(A) + 11t\p(A)e_{\psi}(A)\right]
           {e}^{-(\zeta_A-11e_{\psi}(A)) t\p(A)}. \nonumber
\end{eqnarray}
}  
\\
Term (\ref{16}) is bounded by
\[(\ref{16})\leq \frac{t^k}{k!}{\p(A)}^k(k+1)2n\p(A)
{e}^{-t\p(A)}{e}^{2(k+1)n\p(A)}.\]
To bound term (\ref{17}), we bound the following difference
\[ \left|\frac{\#G_k  k!}{t^k}-1\right|  \leq 
\left|\frac{{\left(t-k(4n)\right)}^k}{t^k}-1\right|  \leq 
 \frac{k\left(k+4n\right)}{t}. \]
Then, we have
\[
(\ref{17})\leq \frac{k\left(k+4n\right)}{t} \frac{{e}^{-t\p(A)}{(t\p(A))}^k}{k!} .\]

Now, we just have to add the five bounds to obtain the theorem with the constant
 $C_{\psi}=1+C_{1}+2C_{p}+8+8$.  Proposition~\ref{18} shows that $C_{1}=5$
and Proposition \ref{mprop} with Theorem \ref{ht} that $C_{p}=116$ .
Then, we prove the theorem with $C_{\psi}=254$.
\end{pf}

\section{Biological applications}\label{BIO}

With the explicit value of the constant $C_{\psi}$ of Theorem \ref{bigt2}, and
more particularly thanks to all the  intermediary bounds given in the proof of this theorem,
we can develop an algorithm to apply this formula to the study of rare words in biological
sequences.
In order to compare different methods, we also compute the bounds corresponding to a $\phi$-mixing,
process for which 
a proof of Poisson approximation is given in \citet{abadi4}. Let us recall the definition of such a
mixing process. 
\begin{defn}
Let $\phi={(\phi(\ell))}_{\ell\geq{0}}$ be a sequence decreasing
to zero. We say that 
 ${(X_m)}_ {m\in{\Zset}}$ is a  $\phi$-mixing process if for all integers $\ell\geq{0}$, the following
holds
\[\sup_{n\in{\Nset},B\in{{\mathcal{F}}_{\{0,.,n\}}},C\in{{\mathcal{F}}_{\{n\geq{0}\}}}}
\frac{|\p(B\cap{T^{-(n+\ell+1)}(C)})-\p(B)\p(C)|}
{\p(B)}=\phi(\ell),\]
where the supremum is taken over the sets $B$ and $C$, such that $\p(B)>0$.
\end{defn}
Note that obviously, $\psi$-mixing implies $\phi$-mixing.
Then, we obtain two new methods for the detection of over- or under-represented words in biological
sequences and we compare them to the Chen-Stein method.

We recall that Markov models are $\psi$-mixing processes and then also $\phi$-mixing
processes. Then, we first need to know the 
functions $\psi$ and $\phi$ for a  Markov model. It turns out that we can use
\[\psi(\ell)=\phi(\ell)=K\nu^{\ell} \mbox{ with $K>0$ and $0<\nu<1$,}\] 
where $K$ and $\nu$ have to be estimated (see \citet{Meyn}). There are several estimations of $K$ and $\nu$. 
We choose $\nu$ equal to the second eigenvalue of the transition matrix of the model
and $K={\left(\inf_{j\in\{1, \dots, {|\mathcal{A}|}^k \}} \mu_j\right)}^{-1}$ where
$|\mathcal{A}|$ is the alphabet size, $k$ the order of the Markov model and $\mu$ the stationary
distribution of the Markov model. 

We recall that we aim at guessing a relevant biological role of a
word in a sequence using its number of occurrences. Thus we compare the number 
of occurrences expected in the Markov chain that models the
sequence and the observed number of occurrences.
It is recommended to choose a degree of significance $s$ 
to quantify this relevance. We fix arbitrarily a degree of significance and we want to calculate
the smallest  
number of occurrences $u$ necessary for $\p(N>u)<s$, where $N$ is the number of occurrences
 of the studied word. 
If the number of occurrences counted in the sequence is larger
than this $u$, we can consider the word to be relevant with a degree of significance $s$. 
We have
\[\p(N>u) \leq \sum_{k=u}^{+\infty}\left(\p_{\mathcal{P}}(N=k)+Error(k)\right)\]
where $\p_{\mathcal{P}}(N=k)$  is the probability under the
Poisson model that $N$ is equal to $k$ and $ Error(k)$ is the error between the exact distribution
and its Poisson approximation, bounded using Theorem {\ref{bigt2}}. 
Then, we search the smallest threshold $u$ such that 
\begin{eqnarray}
\sum_{k=u}^{+\infty}\left(\p_{\mathcal{P}}(N=k)+Error(k)\right)<s. \label{error}
\end{eqnarray}
Then, we have $\p(N>u)<s$
and we consider the word relevant with a degree of significance $s$ if it
appears more than $u$ times in the sequence. 

In order to compare the different methods, we compare
the thresholds that they give.  Obviously, the smaller 
the degree of significance, the more relevant the studied word is.
But for a fixed degree of significance, the best method is the one which 
gives the smallest threshold $u$. Indeed, to give the smallest $u$ is equivalent to give the
smallest error in the tail of the distribution between the exact distribution of the number of
occurrences of word $A$ and the Poisson distribution with parameter $t\p(A)$.

\subsection{Software availability}\label{PANOW}

We developed \texttt{PANOW}, dedicated to the determination of threshold $u$ for given words.
This software is written in ANSI \verb$C++$ and developed on x86
GNU/Linux systems with GCC 3.4, and successfully tested with GCC
latest versions on Sun and Apple Mac OSX systems. It relies 
on \verb$seq++$ library (\citet{seq++}).

Compilation and installation are compliant with the GNU standard procedure.
It is available at \texttt{http://stat.genopole.cnrs.fr/software/panowdir/}. On-line 
documentation is also available.
\texttt{PANOW} is licensed under the GNU General
Public License (\texttt{http://www.gnu.org/licenses/licenses.html}).

\subsection{Comparisons between the three different methods}

\subsubsection{Comparisons using synthetic data.} 

We can compare the mixing methods and the Chen-Stein method through the values of threshold~$u$
obtained with \texttt{PANOW} using \citet{abadi1} in the first case and
\citet{ss2} in the second one.
We recall that the method which gives the smallest threshold $u$ is the best method for a fixed
degree of significance.
Table \ref{t1} offers a good outline of the possibilities and limits of each method. It
displays some results on different words randomly selected (no biological meaning for any of these words).
\begin{table}[htbp]
  \begin{center}
    \caption{Table of thresholds $u$ obtained by the three methods (sequence length $t$  equal to $10^6$).
      {\small
        For each one of the three methods and for each word, we compute the threshold which permits 
        to consider the word as an over-represented word or not, for degree of significance $s$ equal to $0.1$
        or $0.01$. IMP means that the method can not return a result.
      }}
    \label{t1}
    \begin{tabular}{ccccccc}
      &&&&&&\\
      \hline
      & \multicolumn{6}{c}{$t=10^6$}  \\
      \cline{2-7}
      Words & \multicolumn{3}{c}{$s = 0.1$} & \multicolumn{3}{c}{$s = 0.01$} \\
      \cline{2-7}
      & CS & $\phi$ & $\psi$  & CS & $\phi$ & $\psi$  \\
      \hline
      cccg     & IMP  & IMP    & IMP   & IMP  & IMP    & IMP    \\\hline
      aagcgc   & IMP  & $1301$ & $378$ & IMP  & $1304$ & $392$  \\\hline
      cgagcttc & $18$ & $38$   & $18$  & IMP  & $40$   & $22$   \\\hline
      ttgggctg & $14$ & $27$   & $14$  & $18$ & $29$   & $17$   \\\hline
      gtgcggag & $16$ & $32$   & $16$  & $22$ & $34$   & $20$   \\\hline
      agcaaata & $19$ & $39$   & $19$  & IMP  & $41$   & $23$   \\\hline
    \end{tabular}
  \end{center}
\end{table}              
Table \ref{t1} has been obtained with an order one Markov model
using a random transition matrix 
and for a degree of significance of $0.1$ and $0.01$. IMP means that the method can not return a
result. There are several reasons for that and we explain them in
the following paragraph. 
Analysing many results, we notice some differences between the methods.   

Firstly, none of the methods gives us a result in all the cases.
We recall that the Chen Stein method gives a bound ($CS$) using the total variation distance. If the degree
of significance $s$ that we choose is smaller than the bound of Chen-Stein, we never find a threshold
$u$ such that
\[CS + \sum_{k=u}^{+\infty}\p_{\mathcal{P}}(N=k)<s.\]
Then, each time that the given bound is higher than the significance degree, use of the Chen Stein
method is impossible.
Therefore there are many examples that we can not study with
this method. Obviously, it is interesting to have a small degree of significance $s$ and that may be
impossible by this restriction of the Chen-Stein method.
For example, this problem appears for the words
\texttt{aagcgc} and \texttt{cgagcttc} in Table \ref{t1}. For this second word, the Chen-Stein bound
is equal to $0.0107954$. Hence, we can use this method for a significance degree $s$ equal to $0.1$
but not for a significance degree of $0.01$.
The same phenomena appears for the word 
\texttt{agcaaata} (the Chen-Stein bound is equal to $0.0120193$).

The $\phi$- and $\psi$-mixing methods are not based on the total variation distance.
Then, whatever the degree of significance $s$ and if the studied word satisfies the three following
weak properties, we always give a threshold $u$, contrary to the Chen
Stein method. 
In spite of these three conditions, our methods enable us to study a much broader panel of words
than the Chen-Stein method. Indeed, for these two methods, the only problematic cases arise either
when function 
$e_{\psi}$ (see Theorem \ref{bigt2}) is larger 
than $1$ or for a ``high'' parameter of the Poisson distribution (``high'' means larger than $500$)
or when the word periodicity is smaller than half its length (see assumptions in Theorem~\ref{bigt2}: $A
\notin \mathcal{B}_n$).
In fact, the first case does not occur very frequently (in any case in Table~\ref{t1}). The reason
why the function $e_{\psi}$ (or a similar function in the $\phi$-mixing case)
has to be smaller than $1$ is that, for numerical reasons, the error term has to be decreasing with
the number of 
occurrences $k$ and without this condition on $e_{\psi}$ we can not ensure this decrease. We have to
compute error 
terms for a finite number of values of $k$ but in order to reduce the computation time, when error
term becomes smaller than a certain value (we choose $10^{-300}$), we suppose all the following
error terms equals to this value. That is why error term has to be decreasing.
The second problem, a ``high'' parameter of the Poisson distribution, is just a computational
difficulty and once again it does not occur very frequently 
(only for the word \texttt{cccg} in Table \ref{t1} for instance).
We would like to insist on the main advantage of our methods: we can fix any significance degree $s$
and, except in the very rare cases mentioned above, we will find a threshold
$u$, contrary to the Chen-Stein method. 

Also, we can use our methods for any Markov chain order. Indeed,
\texttt{PANOW} runs fast enough contrary to the R program used to compute the Chen-Stein bound of
\citet{ss2}. 
Note that, in program \texttt{PANOW}, we give another method to compute the Chen-Stein
bound (see \citet{abadiT}) and this method gives approximately the same Chen-Stein bound.

The second main observation we can make is that, when it works,
the Chen-Stein method gives either a similar threshold $u$ than the $\psi$-mixing method, or a
smaller one. This means that the $\psi$-mixing method out-performs the Chen-Stein method.

Thirdly we notice that the $\psi$-mixing method is always better than the $\phi$-mixing
one. Obviously, 
this result was expected by the definitions of these mixing processes and also by the theorems because
of the extra factor  
${e}^{-(t-(3k+1)n)\p(A)}$ (see Theorem~\ref{bigt2} and Theorem 2  in \citet{abadi4}). We are
interested by the real impact of this factor on the threshold~$u$: it is significantly better in the case of
a $\psi$-mixing process.

\subsubsection{Biological comparisons.}

Now, we present a few results obtained on real biological examples with order one Markov
models. There are many categories of words 
which have relevant biological functions (promoters, terminators, repeat sequences, chi sites,
uptake sequences, bend sites, signal peptides, binding sites, restriction sites,\~\dots). Some of
them are highly present in the sequence, some others are almost absent. Then, it turns out to be
interesting to
consider the over or the under-representation of words to find words biologically relevant. 

In this section, we test our methods on words already known to be relevant. We focus our study on
Chi sites or uptake sequences. Chi sites of bacterias protect the genome by stopping its
degradation performed by a particular enzyme. The function of this enzyme is to destroy viruses which
could appear into the bacteria. Viruses do not contain Chi sites and then are exterminated. It turns
out that Chi sites are highly present in the bacterial genome.
Uptake sequences are abundant sequence motifs, often located downstream of ORFs, that are used to
facilitate the within-species horizontal transfer of DNA. 

\textit{Example 1}\\
First, we consider the Chi of \textit{Escherichia coli}, \texttt{gctggtgg},
(see Table \ref{coli}), for different degrees
of significance. We use complete sequence of \textit{Escherichia coli K12}
(\citet{ecoliK12}). Sequence length is equal to $4639221$.
\begin{table}[htbp]
  \begin{center}
    \caption{Table of thresholds $u$ obtained by the three methods for the Chi of \textit{Escherichia
        coli}: \texttt{gctggtgg} (sequence length $t$  equal to $4639221$).
      {\small
        For each one of the three methods we compute the threshold which permits 
        to consider the word as an over-represented word or not, for degree of significance $s$. IMP means
        that the method can not return a result. ``counts'' correspond to the number
        of occurrences observed in the sequence.
      }
    }
    \label{coli}
    \begin{tabular}{ccccc}
      &&&&\\
      \hline
      $s$ & Chen-Stein & $\phi$-mixing & $\psi$-mixing & counts\\
      \hline
      $0.1$ & $87$ & $193$ & $83$ & $499$\\
      \hline
      $0.01$ & IMP & $195$ & $92$ & $499$\\
      \hline
      $0.0001$ & IMP & $197$ & $99$ & $499$\\
      \hline
      ${10}^{-239}$& IMP & $549$ & $498$ & $499$\\
      \hline
    \end{tabular}
  \end{center}
\end{table}
We recall that for a fixed significance degree, the smaller the threshold $u$, the best the method is.
Then, we can conclude that the $\psi$-mixing method gives the most interesting results. Chi of
\textit{E. coli} could 
be considered as an over-represented one from 
$99$ occurrences for a significance degree $s$ of $0.0001$. Because Chen-Stein bound is equal to
$0.067726$, Chen-Stein method does not permit to 
conclude for significance degrees of $0.01$ and $0.001$. Moreover, it is well known that
Chi of \textit{E. coli} is a very relevant word in this bacteria. Then, we expect a very small
significance degree for this word. Unfortunately, the minimal significance degree
which could be obtained by Chen-Stein method is, in fact, the Chen-Stein bound: $0.067726$. Our
method allows to obtain very small significance degree and the minimal significance degree for which
Chi of \textit{E. coli} is considered as an over-represented word by the $\psi$-mixing method, is given
at the last line of Table \ref{coli}: it is equal to ${10}^{-239}$.
Note also that the thresholds $u$ increase with the significance degrees $s$. To understand this fact, it
is sufficient to look at inequality (\ref{error}).
But they increase slowly while significance degrees $s$ decreases. It could be surprising but it
is due to the error term which decreases very fast from a certain number of occurrences.

\textit{Example 2}\\
Second, we consider the Chi of \textit{Haemophilus influenzae} and its uptake sequence
 (see Table \ref{haem}), 
for a significance degree $s$ equal to $0.01$. We use complete sequence of \textit{Haemophilus
  influenzae} (\citet{Haeminf}). Sequence length is equal to $1830138$.
\begin{table}[htbp]
  \begin{center}
    \caption{Table of thresholds $u$ obtained by the three methods for the Chi and the uptake sequence
      of \textit{Haemophilus influenzae} (sequence length $t$  equal to $1830138$).
      {\small
        For each one of the three methods and for each word, we compute the threshold which permits 
        to consider the word as an over-represented word or not, for degree of significance equal to
        $0.01$. IMP means that the method can not return a result. ``counts'' correspond to the number
        of occurrences observed in the sequence.
      }
    }
    \label{haem}
    \begin{tabular}{ccccc}
      &&&& \\
      \hline
      Words & Chen-Stein & $\phi$-mixing & $\psi$-mixing & counts\\
      \hline
      gatggtgg (chi)    & $23$ & $36$ & $22$ & $20$  \\
      \hline
      gctggtgg (chi)    & $21$ & $32$ & $20$ & $44$  \\
      \hline
      ggtggtgg (chi)    & $16$ & IMP  & IMP  & $57$  \\
      \hline
      gttggtgg (chi)    & $30$ & $45$ & $26$ & $37$  \\
      \hline
      aagtgcggt (uptake)& $13$ & $17$ & $13$ & $737$ \\
      \hline
    \end{tabular}
  \end{center}
\end{table}
We observe that in all the cases the $\psi$-mixing method is the best one because it gives the
smallest $u$, except for the word
\texttt{ggtggtgg} which has a periodicity less than $\left[\frac{n}{2}\right]$ (and then we can not
study it: see assumptions in Theorem \ref{bigt2}).  We can not assume 
the good significance of the first Chi (\texttt{gatggtgg}) because we count only $20$ occurrences
 in the sequence, 
whereas $23$ occurrences are necessary to consider this word as exceptional. On the other
 hand, the uptake sequence
is very significant (and then very relevant). Indeed, we could fix a significance degree equal to
${10}^{-224}$ and consider it as an over-represented word from $736$ occurrences with the
$\psi$-mixing method. As \texttt{aagtgcggt} is counted $737$ times in the sequence, we obtain the
well-known fact that this word is biologically relevant.

\section{Conclusions and perspectives}
To conclude this paper, we recall the advantages of our new methods. We give an error valid for
all the values $k$ of the random variable $N^t$ corresponding to the number of occurrences of word
$A$ in a sequence of length $t$. Then, we can find a minimal number of occurrences to consider a word as
biologically relevant 
for a very large number of words and for all degrees of significance. That is the main advantage of
our methods on the Chen-Stein one which is based on the total variation 
distance and for which small degrees of significance can not be obtained. Results of our
$\psi$-mixing method and the Chen-Stein method remain similar but our method has
less limitations. Note that our methods provide performing
results for general modelling processes such as Markov chains as
well as every  $\phi$- and $\psi$-mixing processes.

In terms of perspectives, as we expect more significant results, we hope to improve these methods
adapting them 
directly to Markov chains instead of $\psi$- or $\phi$-mixing. 
Moreover, it is well-known that a compound Poisson approximation
is better for self-overlapping words (see \citet{ss1} and \citet{ss2}).
An error term for the compound Poisson
approximation for self-overlapping words can be easily derived
from our results.

\begin{ack}
The authors would like to thank Bernard Prum for his support and his useful comments.
The authors would like to thank Sophie Schbath for her program, Vincent Miele for his very
relevant help in the conception of the software and Catherine Matias for her invaluable advices.
\end{ack}

\bibliographystyle{plainnat}
{
  \footnotesize 
  \bibliography{/home/nvergne/nico/bibliographie/BIBLIO} 

\begin{thebibliography}{33}
\providecommand{\natexlab}[1]{#1}
\providecommand{\url}[1]{\texttt{#1}}
\expandafter\ifx\csname urlstyle\endcsname\relax
  \providecommand{\doi}[1]{doi: #1}\else
  \providecommand{\doi}{doi: \begingroup \urlstyle{rm}\Url}\fi

\bibitem[Abadi(2001{\natexlab{a}})]{abadi2}
M.~Abadi.
\newblock {Exponential approximation for hitting times in mixing processes}.
\newblock \emph{{Mathematical Physics Electronic Journal}}, 7,
  2001{\natexlab{a}}.

\bibitem[Abadi(2004)]{abadi3}
M.~Abadi.
\newblock {Sharp error terms and necessary conditions for exponential hitting
  times in mixing processes}.
\newblock \emph{{Annals of Probability}}, 32:\penalty0 243--264, 2004.

\bibitem[Abadi(2001{\natexlab{b}})]{abadiT}
M.~Abadi.
\newblock \emph{{Instantes de ocorr\^{e}ncia de eventos raros em processos
  misturadores}}.
\newblock PhD thesis, Universidade de S\~{a}o paulo, 2001{\natexlab{b}}.
\newblock available at \texttt{http://www.ime.unicamp.br/\~{}miguel}.

\bibitem[Abadi and Vergne(2006{\natexlab{a}})]{abadi1}
M.~Abadi and N.~Vergne.
\newblock {Sharp error terms for return time statistics under mixing
  conditions}.
\newblock Submitted, 2006{\natexlab{a}}.

\bibitem[Abadi and Vergne(2006{\natexlab{b}})]{abadi4}
M.~Abadi and N.~Vergne.
\newblock {Sharp error terms for {P}oisson statistics under mixing conditions:
  {A} new approach}.
\newblock Submitted, 2006{\natexlab{b}}.

\bibitem[Almagor(1983)]{MC}
H.~Almagor.
\newblock {A} {M}arkov analysis of {DNA} sequences.
\newblock \emph{J.Theor. Biol.}, 104:\penalty0 633--645, 1983.

\bibitem[Arratia et~al.(1989)Arratia, Goldstein, and Gordon]{Ar1}
R.~Arratia, L.~Goldstein, and L.~Gordon.
\newblock {T}wo moments suffice for {P}oisson approximations: the
  {C}hen-{S}tein method.
\newblock \emph{Ann. Prob.}, 17:\penalty0 9--25, 1989.

\bibitem[Arratia et~al.(1990)Arratia, Goldstein, and Gordon]{Ar2}
R.~Arratia, L.~Goldstein, and L.~Gordon.
\newblock {P}oisson approximation and the {C}hen-{S}tein method.
\newblock \emph{Statist. Sci.}, 5:\penalty0 403--434, 1990.

\bibitem[Barbour et~al.(1992)Barbour, Chen, and Loh]{Bar}
A.D. Barbour, L.H.Y. Chen, and W.L. Loh.
\newblock {C}ompound {P}oisson approximation for nonnegative random variables
  via {S}tein's method.
\newblock \emph{Ann. Prob.}, 20:\penalty0 1843--1866, 1992.

\bibitem[Blaisdell(1985)]{MCA}
B.E. Blaisdell.
\newblock {M}arkov chain analysis finds a significant influence of neighboring
  bases on the occurrence of a base in eucaryotic nuclear {DNA} sequences both
  protein-coding and noncoding.
\newblock \emph{J. Mol. Evol.}, 21:\penalty0 278--288, 1985.

\bibitem[Blattner et~al.(1997)Blattner, Plunkett, Bloch, Perna, Burland, Riley,
  Collado-Vides, Glasner, Rode, G.F., Gregor, N.W., Kirkpatrick, Goeden, Rose,
  Mau, and Shao]{ecoliK12}
F.R. Blattner, G.3rd Plunkett, C.A. Bloch, N.T. Perna, V.~Burland, M.~Riley,
  J.~Collado-Vides, J.D. Glasner, C.K. Rode, Mayhew G.F., J.~Gregor, Davis
  N.W., H.A. Kirkpatrick, M.A. Goeden, D.J. Rose, B.~Mau, and Y.~Shao.
\newblock The complete genome sequence of escherichia coli k-12.
\newblock \emph{Science}, 277:\penalty0 1453--74, 1997.

\bibitem[Chen(1975)]{Chen}
L.H.Y. Chen.
\newblock {P}oisson approximation for dependant trials.
\newblock \emph{Ann. Prob.}, 3:\penalty0 534--545, 1975.

\bibitem[Douglass(1996)]{MVT}
S.A. Douglass.
\newblock \emph{Introduction to Mathematical Analysis}, chapter~8.
\newblock Addison-Wesley, Boston, 1996.

\bibitem[El~Karoui et~al.(1999)El~Karoui, Biaudet, Schbath, and Gruss]{ElK}
M.~El~Karoui, V.~Biaudet, S.~Schbath, and A.~Gruss.
\newblock {C}haracteristics of {C}hi distribution on different bacterial
  genomes.
\newblock \emph{Res. Microbiol.}, 150:\penalty0 579--587, 1999.

\bibitem[Fleischmann et~al.(1995)Fleischmann, Adams, White, and
  Clayton]{Haeminf}
R.D. Fleischmann, M.D Adams, O.~White, and R.A. Clayton.
\newblock Whole-genome random sequencing and assembly of haemophilus influenzae
  rd.
\newblock \emph{{Science}}, 269:\penalty0 496--512, 1995.

\bibitem[Gelfand et~al.(1992)Gelfand, Kozhukhin, and P.A.]{Extend}
M.S. Gelfand, C.G. Kozhukhin, and Pevzner P.A.
\newblock {E}xtendable words in nucleotide sequences.
\newblock \emph{Bioinformatics}, 8:\penalty0 129--135, 1992.

\bibitem[Godbole(1991)]{Poisson}
A.P. Godbole.
\newblock {P}oisson approximations for runs and patterns of rare events.
\newblock \emph{Adv. Appl. Prob.}, 23:\penalty0 851--865, 1991.

\bibitem[Karlin et~al.(1992)Karlin, Burge, and Campbell]{Karlin2}
S.~Karlin, C.~Burge, and A.M. Campbell.
\newblock Statistical analyses of counts and distributions of restriction sites
  in dna sequences.
\newblock \emph{{Nucl. Acids Res.}}, 20:\penalty0 1363--1370, 1992.

\bibitem[Meyn and Tweedie(1993)]{Meyn}
S.P. Meyn and R.L. Tweedie.
\newblock \emph{{M}arkov {C}hains and {S}tochastic {S}tability}.
\newblock Springer-Verlag, Heidelberg, 1993.

\bibitem[Miele et~al.(2005)Miele, Bourguignon, Robelin, Nuel, and
  Richard]{seq++}
V.~Miele, P.Y. Bourguignon, D.~Robelin, G.~Nuel, and H.~Richard.
\newblock {seq++~: analyzing biological sequences with a range of
  Markov-related models}.
\newblock \emph{Bioinformatics}, 21:\penalty0 2783--2784, 2005.

\bibitem[Nicod\`{e}me et~al.(2002)Nicod\`{e}me, Doerks, and Vingron]{Nicod}
P.~Nicod\`{e}me, T.~Doerks, and M.~Vingron.
\newblock Proteome analysis based on motif statistics.
\newblock \emph{{Bioinformatics}}, 18\penalty0 (Suppl. 2):\penalty0 5161--5171,
  2002.

\bibitem[Nuel(2004)]{GregSpatt}
G.~Nuel.
\newblock {LD-SPatt: Large Deviations Statistics for Patterns on Markov
  chains}.
\newblock \emph{Comp. Biol.}, 11:\penalty0 1023--1033, 2004.

\bibitem[Phillips et~al.(1987)Phillips, Arnold, and Ivarie]{hexa}
G.J. Phillips, J.~Arnold, and R.~Ivarie.
\newblock {T}he effect of codon usage on the oligonucleotide composition of the
  e. coli genome and identification of over- and underrepresented sequences by
  {M}arkov chain analysis.
\newblock \emph{Nucl. Acids Res.}, 15:\penalty0 2627--2638, 1987.

\bibitem[Prum et~al.(1995)Prum, Rodolphe, and de~Turckheim]{Bernard2}
B.~Prum, F.~Rodolphe, and E.~de~Turckheim.
\newblock Finding words with unexpected frequencies in {DNA} sequences.
\newblock \emph{{J. R. Statis. Soc. B}}, 11:\penalty0 190--192, 1995.

\bibitem[R\'{e}gnier(2000)]{MR}
M.~R\'{e}gnier.
\newblock {A} unified approach to word occurrence probabilities.
\newblock \emph{Discr. Appl. Math.}, 104:\penalty0 259--280, 2000.

\bibitem[Reinert and Schbath(1998)]{ss2}
G.~Reinert and S.~Schbath.
\newblock {C}ompound {P}oisson and {P}oisson process approximations for
  occurrences of multiple words in {M}arkov chains.
\newblock \emph{J. Comput. Biol.}, 5:\penalty0 223--253, 1998.

\bibitem[Reinert et~al.(2000)Reinert, Schbath, and Waterman]{ss1}
G.~Reinert, S.~Schbath, and M.S. Waterman.
\newblock {P}robabilistic and {S}tatistical {P}roperties of {W}ords: {A}n
  {O}verview.
\newblock \emph{J. Comput. Biol.}, 7, 2000.

\bibitem[Robin and Daudin(1999)]{jjd}
S.~Robin and J.J. Daudin.
\newblock {E}xact distribution of word occurrences in a random sequence of
  letters.
\newblock \emph{J. Appl. Prob.}, 36, 1999.

\bibitem[Smith et~al.(1981)Smith, Kunes, Schultz, Taylor, and Triman]{ecoli}
G.R. Smith, S.M. Kunes, D.W. Schultz, A.~Taylor, and K.L. Triman.
\newblock {S}tructure of chi hotspots of generalized recombination.
\newblock \emph{Cell}, 24:\penalty0 429--36, 1981.

\bibitem[Smith et~al.(1999)Smith, Gwinn, and Salzberg]{uptake}
H.O. Smith, M.L Gwinn, and S.L. Salzberg.
\newblock {DNA} uptake signal sequences in naturally transformable bacteria.
\newblock \emph{Res. Microbiol.}, 150:\penalty0 603--616, 1999.

\bibitem[Stein(1972)]{Stein}
C.~Stein.
\newblock {A} bound for the error in the normal approximation to the
  distribution of a sum of dependent random variables.
\newblock \emph{{Proc. Sixth Berkeley Symp. Math. Statist. Probab.}},
  2:\penalty0 583--602, 1972.
\newblock University of California Press.

\bibitem[van Helden et~al.(1998)van Helden, Andr\'{e}, and
  Collado-Vides]{Van_H2}
J.~van Helden, B.~Andr\'{e}, and J.~Collado-Vides.
\newblock Extracting regulatory sites from the upstream region of yeast genes
  by computational analysis of oligonucleotide frequencies.
\newblock \emph{{J. Mol. Biol.}}, 281:\penalty0 872--842, 1998.

\bibitem[van Helden et~al.(2000)van Helden, del Olmo, and
  P\'{e}rez-Ort\'{i}n]{Van_H}
J.~van Helden, M.~del Olmo, and J.E. P\'{e}rez-Ort\'{i}n.
\newblock Statistical analysis of yeast genomic downstream sequences reveals
  putative polyadenylation signals.
\newblock \emph{{Nucl. Acids Res.}}, 28:\penalty0 1000--1010, 2000.

\end{thebibliography}
}
\end{document}